\documentclass[11pt,a4paper]{article}
\usepackage{amsmath, amssymb, amsthm}
\usepackage[left=1in, right=1in, top=1in, bottom=1in, includefoot, headheight=14pt]{geometry}
\usepackage[colorlinks]{hyperref}
\usepackage{stackrel,setspace}
\usepackage{tikz}
\newcommand\beq{\begin{equation}}
\newcommand\eeq{\end{equation}}

\newcommand{\bkE}{\mathbb{E}}

\newcommand{\p}{{\mathbb{P}}}

\DeclareMathOperator{\sign}{sign}

\newcommand{\tMonm}{\widetilde{\mathcal{F}}}
\newcommand{\tMon}{\widetilde{{\rm Mon}}}
\numberwithin{equation}{section}
\newcommand{\tq}{\,:\,}

\newcommand{\I}{{\mathbf 1}}

\newcommand{\E}{{\mathbb{E}}}

\newcommand{\bb}{\begin{eqnarray*}}
\newcommand{\ee}{\end{eqnarray*}}
\newcommand{\bbb}{\begin{eqnarray}}
\newcommand{\eee}{\end{eqnarray}}
\newcommand{\F}{{\mathcal{F}}}

\theoremstyle{plain}

\newtheorem{theorem}{Theorem}[section]
\newtheorem{Theorem}{Theorem}[section]

\newtheorem{Proposition}[theorem]{Proposition}

\newtheorem{lma}[theorem]{Lemma}
\newtheorem{Definition}[theorem]{Definition}

\theoremstyle{definition}
\newtheorem{definition}[theorem]{Definition}
\newtheorem{remark}[theorem]{Remark}

\allowdisplaybreaks[3]

\begin{document}

\begin{center}
{\Large \textbf{A deviation bound for $\alpha$-dependent sequences with applications to intermittent maps}} \medskip

\bigskip
J. Dedecker$^a$ and F. Merlev\`{e}de$^b$
\bigskip

\end{center}

$^a$ Universit\'e Paris Descartes, Sorbonne Paris Cit\'e,  Laboratoire MAP5 (UMR 8145).\\
Email: jerome.dedecker@parisdescartes.fr \bigskip

$^b$ Universit\'{e} Paris Est, LAMA (UMR 8050), UPEM, CNRS, UPEC.\\
Email: florence.merlevede@u-pem.fr \bigskip

\bigskip



\begin{center}
{\bf Abstract}
\end{center}

We prove a deviation bound for the maximum of partial 
sums of functions of  $\alpha$-dependent sequences as defined in \cite{DGM}.
As a consequence, we extend the Rosenthal 
inequality of Rio \cite{Ri} for $\alpha$-mixing sequences in the sense of Rosenblatt \cite{R} to the larger class of $\alpha$-dependent 
sequences. Starting from the deviation inequality, we obtain upper bounds 
for large deviations and an H\"olderian invariance principle for the 
Donsker line. We illustrate our results through the example of
intermittent maps of the interval, which are not $\alpha$-mixing
in the sense of Rosenblatt.

\section{Introduction} 
For stationary $\alpha$-mixing sequences in the sense of Rosenblatt (see \cite{R}) a 
Fuk-Nagaev type inequality has been proved by Rio (see Theorem 6.2 in \cite{Ri}). This deviation 
inequality is very powerful and enables one to prove optimal results for the deviation of partial sums 
and a sharp Rosenthal type inequality   (see Corollary 6.3
in \cite{Ri}). The proof uses a 
blocking technique and the coupling property of $\alpha$-mixing coefficients. 
 
Rio's inequality has been extended to a larger class of dependent sequences
in \cite{DP}. In that paper, the authors  noticed that one can use a dependency coefficient $\tau(n)$ 
whose definition is perfectly adapted to the coupling property, in the spirit
of R\"uschendorf \cite{Ru}. The Fuk-Nagaev inequality
for $\tau$-dependent sequences applies to many non-mixing sequences in the sense of Rosenblatt, 
such as contracting Markov chains or causal functions of infinite sequences of independent and
identically distributed (i.i.d.) random variables. 

However, although the property of $\tau$-dependency  is much less
restrictive than $\alpha$-mixing, it appears   to  be not
well adapted to most of the usual  dynamical systems. The main reason is that, 
to prove the Fuk-Nagaev inequality {\em via} blocking + coupling, 
one needs to control the dependency between the whole past 
and the whole future of the sequence. To the best of our
 knowledge, this
can be done only for a very restricted class of dynamical systems 
(see  Example 4 of Section 7.2  in \cite{DP2}). 

In the present paper, we prove a deviation bound for the maximum of partial 
sums of functions of stationary $\alpha$-dependent sequences as defined in \cite{DGM}.
More precisely, if ${\bf X}=(X_i)_{i \in {\mathbb Z}}$ is a strictly stationary sequence of
real-valued random variables, our deviation inequality (see Proposition \ref{deviationRosenthal}) is expressed in terms of
 a coefficient $\alpha_{2, {\bf X}}(n)$ which only controls the dependency between the past
 (before time 0) 
 and any variable of the form 
 $
 {\mathbf 1}_{X_i \leq t, X_j \leq s}
 $
 where $i, j \geq n$. Note that this coefficient can be exactly computed 
 for the  Markov chains associated to the  intermittent maps 
 introduced
  in \cite{DGM} (see Subsection \ref{subsec:int}).  We shall also describe precisely 
 the class of observables $f(X_i)$ to which our results apply
 (in particular, this class  contains the functions $f$ 
 which are  piecewise monotonic on open intervals with a finite number of branches, under an appropriate control of the 
 tail of $f(X_0)$).  
 
The proof of our deviation inequality still uses a blocking argument, but the coupling part is now replaced by a martingale approximation
followed by an application of the Rosenthal-type inequality 
proved in \cite{MP}.  The resulting  inequality is not of the same kind as the usual
 Fuk-Nagaev inequality, but it seems to perform as well in most of
 the applications. For instance, it provides a full  extension of the Rosenthal inequality of Rio \cite{Ri} to the larger class of 
 $\alpha$-dependent sequences (see our Theorem \ref{ineRosalphafaible}). We also use it  to obtain upper bounds 
for large deviations and an H\"olderian invariance principle for the 
Donsker line. Concerning the H\"olderian invariance principle, we follow the approach of Giraudo
\cite{Gir}, who recently obtained very precise 
results for mixing sequences ($\alpha$-mixing  in the sense of Rosenblatt, $\tau$-dependency and $\rho$-mixing)
by applying sharp deviation inequalities for the maximum 
of partial sums. 

The paper is organized as follows. 
In Section \ref{DefNot},  we give the notations and definitions which will be used all
along the paper.  
In Section \ref{Sec:new}, we present the main 
consequences of our deviation inequalities: moment 
bounds in Subsection \ref{Sec:MR}, large deviation bounds
in Subsection \ref{Sec:LD}, H\"olderian invariance 
principles in Subsection \ref{Sec:HIP}. The application 
of these results to intermittent maps are given in Section
\ref{subsec:int}. In Section \ref{sectiondeviation} our 
main deviation inequality is stated and proved.
Finally, the proofs of the results of Section \ref{Sec:new}  are gathered in Section \ref{Sec:6}.

\section{Definitions and notations} \label{DefNot}
\setcounter{equation}{0}


Let us start with the  notation
 $a_n(x) \ll b_n(x)$, which means that there exists a positive  constant $C$ not
depending on $n$ nor $x$ such that  $a_n(x) \leq  Cb_n(x)$, for all positive integers  $n$ and all real $x$.
\subsection{Stationary sequences and dependency coefficients}

Let  $(\Omega ,\mathcal{A}, \p)$ be a probability space, and 
$T :\Omega \mapsto \Omega $ be a bijective bi-measurable
transformation preserving the probability ${\p}$. Let
${\mathcal F}_0$ be a sub-$\sigma$-algebra of $\mathcal{A}$
satisfying ${\mathcal F}_0 \subseteq T^{-1}({\mathcal
F}_0)$. 

Let $Y_0$ be an ${\mathcal F}_0$-measurable and integrable real-valued random variable
with distribution $P_{Y_0}$. 
Define the stationary sequence 
 ${\bf Y}=(Y_i)_{i \in {\mathbb Z}}$ by $Y_i=Y_0 \circ T^i$.
 

\medskip


Let us now  define the dependency coefficients of the sequence 
$(Y_i)_{i \in {\mathbb Z}}$ as in \cite{DGM}. These coefficients
are less restrictive than the usual mixing coefficients of Rosenblatt \cite{R}. 

\medskip

\begin{definition}
For any integrable random variable $Z$, let 
$Z^{(0)}=Z- \E(Z)$.
For any random variable $V=(V_1, \cdots, V_k)$ with values in
${\mathbb R}^k$ and any $\sigma$-algebra ${\mathcal F}$, let
\[
\alpha({\mathcal F}, V)= \sup_{(x_1, \ldots , x_k) \in {\mathbb R}^k}
\left \| \E \left(\prod_{j=1}^k (\I_{V_j \leq x_j})^{(0)} \Big | {\mathcal F} \right)
-\E \left(\prod_{j=1}^k (\I_{V_j \leq x_j})^{(0)}  \right) \right\|_1.
\]
For the stationary sequence ${\bf Y}=(Y_i)_{i \in {\mathbb Z}}$, let 
\begin{equation}
\label{defalpha} \alpha_{k, {\bf Y}}(0) = 1/2 \, \text{ and } \, \alpha_{k, {\bf Y}}(n) = \max_{1 \leq l \leq
k} \ \sup_{ n\leq i_1\leq \ldots \leq i_l} \alpha({\mathcal F}_0,
(Y_{i_1}, \ldots, Y_{i_l}))  \, \text{ for $n >0$.}
\end{equation}
Note that $\alpha_{1, {\bf Y}}(n)$ is then simply given by 
\begin{equation}
\alpha_{1, {\bf Y}}(n)= \sup_{x \in {\mathbb R}}\left \|\E\left (\I_{Y_n \leq x}|{\mathcal F}_0\right ) - F(x)
\right \|_1 \, ,
\end{equation} 
where $F$ is the distribution function of $P_{Y_0}$. 
\end{definition}

\medskip

All the results of the paper  involve only the coefficients 
$\alpha_{1, {\bf X}}(n)$ and 
$\alpha_{2, {\bf X}}(n)$.

\subsection{Quantile functions and piecewise monotonic 
spaces}

In this subsection, we describe the functions spaces to which our 
results apply. Contrary to the usual  mixing case, any function of a
stationary  
$\alpha$-dependent  sequence
 ${\bf Y}=(Y_i)_{i \in {\mathbb Z}}$ is not necessarily $\alpha$-dependent (meaning that its dependency 
coefficients do no necessarily tend to zero). Hence, we need to impose some constraints on the observables. 

The first thing to notice is that, if $f$ is monotonic on some open interval and $0$  elsewhere, and if 
${\bf X}=(f(Y_i))_{i \in {\mathbb Z}}$, then for any positive integer $k$,
$$
\alpha_{k, {\bf X}}(n) \leq 2^k \alpha_{k, {\bf Y}}(n) \, .
$$
As a consequence, if one can prove a deviation inequality for $\sum_{k=1}^n Y_i$ with an upper bound involving the coefficients $( \alpha_{k, {\bf Y}}(n))_{n \geq 0}$ then it also holds for $
\sum_{k=1}^n f(Y_i)
$, 
where $f$ is  monotonic on a single interval. In this case, the deviation inequality can be  extended  by linearity to convex combinations of such functions.  Such classes are
described in Definition \ref{defclosedenv}  below.

First, we need a uniform control on the tail of our test functions 
by a given tail function $H$.

\begin{Definition}
A function $H$ from ${\mathbb R}^+$ to $[0,1]$ is 
a tail function if it  is non-increasing, right-continuous and 
converges to zero at infinity. The quantile function 
$Q=H^{-1}$ 
is the  generalized inverse of the tail function $H$:
for $u \in [0,1]$,
$
Q(u) = \inf \left \{  t \geq 0 : H(t) \leq
u\right \} $ (with the convention that 
$\inf\{\emptyset\}= \infty$). For $p \geq 1$,  we say that 
$Q$ belongs to ${\mathbb L}^p$ if  $\int_0^1 Q^p(u) du < \infty$. 
\end{Definition}


The function spaces are then defined as follows: 

\begin{Definition}
\label{defclosedenv} 
If $\mu$ is a probability measure on ${\mathbb R}$ 
and $Q=H^{-1}$ is an integrable quantile function,
let $\tMon(Q, \mu)$ be the set of functions $g$ which are
monotonic on some open interval of ${\mathbb R}$ and null
elsewhere and such that $\mu(|g|>t) \leq H(t)$
for any $t \in {\mathbb R}^+$. Let $\tMonm(Q,
\mu)$ be the closure in ${\mathbb L}^1(\mu)$ of the set of
functions which can be written as $\sum_{\ell=1}^{L} a_\ell
f_\ell$, where $\sum_{\ell=1}^{L} |a_\ell| \leq 1$ and $f_\ell$
belongs to $\tMon(Q, \mu)$.
\end{Definition}

Note that a function belonging to $\tMonm(Q, \mu)$
 is allowed to blow up at an infinite number of
points.  Note also that any function
$f$
with bounded variation (BV) such that
$|f| \leq M_1$
and
$\|
df
\| \leq 
M_2
$
belongs to the class $\tMonm(Q, \mu)$
for any $\mu$
and the quantile function
$Q \equiv M_1+2M_2$
(here $\|
df
\|
$
denotes the variation norm of the signed measure
$df$).  
Moreover,  if a
function
$f$
is piecewise monotonic with
$N$
branches,  then it belongs to
$\tMonm(Q, \mu)$ for the quantile function 
$Q=H^{-1}$ and 
$H(t)= \mu( |f|> t/N)$. 
   Finally,  let  us  emphasize  that  there  is  no  requirement  on  the  modulus  of
continuity for functions in $\tMonm(Q, \mu)$.

\section{Main results for $\alpha$-dependent sequences}
\label{Sec:new}

In Proposition \ref{deviationRosenthal} given in Section \ref{sectiondeviation}, we give a general deviation inequality for $\alpha$-dependent sequences. In this section, we present some 
striking applications of this inequality. 

We shall use the following notations:  for $u \in [0,1]$ and $k \in {\mathbb N}^*$, let 
\begin{equation}\label{invalpha}
 \alpha_{k,{\bf Y}}^{-1} (u)=\min\{q\in {\mathbb N} \tq \alpha_{k,{\bf Y}}(q)\leq u\}
 = \sum_{n=0}^\infty {\bf 1}_{u < \alpha_{k,{\bf Y}}(n)}
 \, .
\end{equation}
Note that $\alpha_{1,{\bf Y}}(n) \leq 
\alpha_{2,{\bf Y}}(n)$, and consequently 
$\alpha_{1,{\bf Y}}^{-1} \leq \alpha_{2,{\bf Y}}^{-1}$. 

\subsection{Moment  inequalities}
\label{Sec:MR}


Our first result is the following Rosenthal-type  inequality for the maximum of  partial sums of $\alpha$-dependent sequences for all powers $p \geq 2$. 

\begin{Theorem} \label{ineRosalphafaible}
Let $p \geq 2$ and 
let $Q$ be a  quantile function in ${\mathbb L}^p$. Let ${\bf Y}=(Y_i)_{i \in {\mathbb Z}}$,
where $Y_i= Y_0 \circ T^i$.
Denote by $P_{Y_0}$  the distribution of $Y_0$. Let  $X_i = f(Y_i) - \E ( f(Y_i))$, where  $f$ belongs to $\tMonm(Q, P_{Y_0})$ and let $S_n =\sum_{k=1}^n X_k$. 
Then
\begin{equation}\label{RB}
\left \Vert
 \max_{ 1 \leq k \leq n}  |
  S_k | 
\right \Vert_p^p \ll n^{p/2} 
\left ( \int_0^1 (\alpha_{1,{\bf Y}}^{-1} (u) \wedge n)   Q^{2}(u) du \right )^{p/2} + n \int_0^1 (\alpha_{2,{\bf Y}}^{-1} (u) \wedge n)^{p-1}   Q^{p}(u) du \, .
\end{equation}
\end{Theorem}

\begin{remark}  Note  that 
Inequality \eqref{RB} can be written as follows:
\begin{equation}\label{RBvar}
\left \Vert
 \max_{ 1 \leq k \leq n}  |
  S_k | 
\right \Vert_p^p \ll n^{p/2} 
\left ( \sum_{k=0}^{n-1} \int_0^{\alpha_{1,
{\bf Y}} (k)}Q^2 (u) du \right )^{p/2} + n \sum_{k=0}^{n-1}
(k+1)^{p-2} \int_0^{\alpha_{2,
{\bf Y}} (k)}   Q^{p}(u) du \, .
\end{equation}
\end{remark}

\begin{remark}
Inequality \eqref{RB} is an extension of
Rio's  inequality 
for $\alpha$-mixing sequences  (see  Theorem 6.3 
in  \cite{Ri}) 
to the larger class of  $\alpha$-dependent sequences 
as defined in  \eqref{defalpha} (with $k=2$ for the index of the dependency). 
 Note that Rio's inequality cannot be applied 
to the class of GPM maps described in Subsection
\ref{subsec:int}, because the associated Markov chains of such  maps are not $\alpha$-mixing  in the sense of 
Rosenblatt. Note also that Inequality \eqref{RBvar} implies in particular that if $p\geq 2$ and \[ \sum_{k \geq 0}
(k+1)^{p/2-1} \int_0^{\alpha_{2,
{\bf Y}} (k)}   Q^{p}(u) du < \infty, \quad  
 \text{then} \quad  \left \Vert
 \max_{ 1 \leq k \leq n}  |
  S_k | 
\right \Vert_p^p  \ll n^{p/2} \, .
  \]
    We refer to Section 6.4  in Rio \cite{Riobis} for other possible consequences of Inequality \eqref{RB} (see in particular Corollary 6.1 
  in Rio \cite{Riobis}).
\end{remark}


\subsection{Large deviation inequalities}
\label{Sec:LD}

In this section, we give some upper bounds for the quantity 
\[
 \p \left(\frac 1 n  \max_{ 1 \leq k \leq n}  |
  S_k | \geq x \right) \, .
\]
We shall use the notation
\begin{equation} \label{defofRuinfinity}
  R(u)=\alpha_{2,{\bf Y}}^{-1} (u) Q(u)\, ,
 \quad 
  \text{for $ u  \in [0,1]$}.
\end{equation}

\begin{Theorem}\label{Prop:LD}
Let $Q$ be a quantile function in ${\mathbb L}^1$, and let
 $Y_i$, $f$, $X_i$ and $S_n$ be as in Theorem \ref{ineRosalphafaible}. 
 \begin{enumerate}
\item Assume that, for $p\geq 2$,
\begin{equation}\label{WM}
 \sup_{x>0}x^{p-1} \int_0^1 Q(u) {\bf 1}_{R(u)>x} du  < \infty \, .
\end{equation}
Then, for $p>2$, any  $a \in (p-1, p)$ and any $x>0$,
\begin{equation}\label{WB}
\ \p \left(\frac 1 n  \max_{ 1 \leq k \leq n}  |
  S_k | \geq x \right)
 \ll \frac 1 {n^a x^{2a}} + \frac 1 {n^{p-1} x^p} \, .
\end{equation}
For $p=2$, any  $a \in (1, 2)$, any $c \in (0,1)$ and any $x>0$,
\begin{equation}\label{WB2}
\quad \p \left(\frac 1 n  \max_{ 1 \leq k \leq n}  |
  S_k | \geq x \right)
 \ll \frac 1 {n^{ac} x^{a(1+c)}} + \frac 1 {n x^2} \, .
\end{equation}
\item Assume that, for $p\geq 2$, 
\begin{equation}\label{SM}
 \int_0^1 (\alpha_{2,{\bf Y}}^{-1} (u))^{p-1}   Q^{p}(u) du 
 < \infty \, .
\end{equation}
Then, for any $a \in (p-1, p)$ and any $x>0$, 
\begin{equation}\label{SB}
\sum_{n >0} n^{p-2}\p \left(\frac 1 n  \max_{ 1 \leq k \leq n}  |
  S_k | \geq x \right)
 \ll \frac 1 { x^{2a}} + \frac 1 { x^p} \, .
\end{equation}
\end{enumerate}
\end{Theorem}


\begin{remark}\label{Rem:p<2}
We consider here the case where   $p \in (1,2)$.
Let $Q$ be a quantile function in ${\mathbb L}^1$, and let
 $Y_i$, $f$, $X_i$ and $S_n$ be as in Theorem \ref{ineRosalphafaible}.
  Following the proof of Theorem \ref{Prop:LD}
and using Proposition 1 in \cite{DM} instead of Inequality \eqref{eqmainmax}, one can prove that:
\begin{enumerate}
\item If \eqref{WM} holds, then for  any $x>0$, 
\begin{equation}\label{WBeasy}
\p \left(\frac 1 n  \max_{ 1 \leq k \leq n}  |
  S_k | \geq x \right)
 \ll  \frac 1 {n^{p-1} x^p} \, .
\end{equation}
\item 
If \eqref{SM} holds, then for  any $x>0$, 
\begin{equation}\label{SBeasy}
\sum_{n >0} n^{p-2}\p \left(\frac 1 n  \max_{ 1 \leq k \leq n}  |
  S_k | \geq x \right)
 \ll  \frac 1 { x^p} \, .
\end{equation}
\end{enumerate}
We refer to Subsection 4.2 in \cite{DGM} to see how to apply
 Proposition 1 in  \cite{DM} to $\alpha$-dependent sequences.
Note that in  the case  $p \in (1,2)$, the conditions \eqref{WM}
and \eqref{SM} can be slightly weakened by replacing 
$\alpha_{2,{\bf Y}}^{-1} (u)$ by $\alpha_{1,{\bf Y}}^{-1} (u)$
in the definition of $R(u)$. 
\end{remark}

\begin{remark}
Since $x^{p-1} {\bf 1}_{R(u)>x} \leq (R(u))^{p-1}$, we see
 that  the condition \eqref{SM} is stronger than 
\eqref{WM}. Note also that  \eqref{SM} is equivalent to
$$
\sum_{k=0}^{\infty}
(k+1)^{p-2} \int_0^{\alpha_{2,
{\bf Y}} (k)}   Q^{p}(u) du < \infty \, .
$$
From \eqref{WB},  \eqref{WB2} and \eqref{WBeasy}, it follows that,  for any $x>0$
and any $p >1$, 
$$\p \left (  \max_{ 1 \leq k \leq n}  |
  S_k | \geq nx \right ) =O\left(\frac {1}{n^{p-1}}\right)\, .
  $$
   From \eqref{SB} and \eqref{SBeasy}, it follows 
  that, for any $x>0$ and any $p>1$,
\begin{equation}\label{petito}
\p \left(  \max_{ 1 \leq k \leq n}  |
  S_k | \geq nx \right ) =o\left(\frac {1}{n^{p-1}}\right)\, . 
 \end{equation}
\end{remark}

\begin{remark}
 Let $m$ be a non-negative integer. 
As usual, the stationary sequence $\bf X$ is 
$m$-dependent if 
$\sigma(X_i, i \leq 0)$ is independent of 
$\sigma(X_i, i \geq m+1)$,  and $m=0$ 
corresponds to  the case of i.i.d  random variables. 
If $\bf X$ is 
a stationary $m$-dependent sequence of centered random variables,
we infer from Theorem \ref{Prop:LD}  (resp.  
Remark \ref{Rem:p<2}) that \eqref{WB}-\eqref{WB2}
 (resp. \eqref{WBeasy}) holds
 for $p\geq 2$ (resp. for $p \in (1,2)$) as soon as 
\begin{equation} \label{mdep}
\sup_{x>0}x^{p-1} {\mathbb E}(|X_0|{\bf 1}_{|X_0|>x}) <\infty \, .
\end{equation} 
 Since $p>1$, It is easy to see
that \eqref{mdep} is equivalent to 
$$
\sup_{x>0}x^{p} {\mathbb P}(|X_0|>x) <\infty \, ,
$$
meaning that $X_0$ has a weak moment of order $p$.
In the same way, \eqref{SB}  (resp. \eqref{SBeasy}) holds
for $p \geq 2$ (resp. for $p \in (1,2)$) as soon as 
 ${\mathbb E}(|X_0|^p)< \infty$. In particular,
  if ${\mathbb E}(|X_0|^p)< \infty$ for $p>1$, then \eqref{petito} holds. Now,  according to 
 Proposition 2.6 in \cite{LV}, the estimate \eqref{petito}
  cannot be 
 essentially improved in the i.i.d. case. 
\end{remark}

\subsection{H\"olderian invariance principles}\label{Sec:HIP}

Let  $Y_i$, $f$, $X_i$ and $S_n$ be as in Theorem \ref{ineRosalphafaible}, and define
$$
W_n(t)= \frac 1 {\sqrt{n}} S_{[nt]} +\frac{(nt-[nt])}{\sqrt n} 
X_{[nt]+1} \, .
$$
From \cite{DGM} we know that, if
\begin{equation}\label{DMR}
 \int_0^1 \alpha_{2,{\bf Y}}^{-1} (u)   Q^{2}(u) du 
 < \infty \, ,
\end{equation}
then the process
$
\left \{ W_n(t), t \in [0,1]\right \}
$
converges in distribution in the space $C([0,1], \|\cdot \|_\infty)$ of continuous bounded function on $[0,1]$ to $\sigma W$, where $W$ is a standard 
Brownian motion and 
\begin{equation}\label{sigma}
\sigma^2= {\mathrm{Var}}(X_0)
+
2 \sum_{k>0} {\mathrm{Cov}}(X_0, X_k)
\, .
\end{equation}

For $\beta \in (0,1]$,  let ${\mathcal H}_\beta([0,1])$ be the 
Banach space of $\beta$-H\"older functions 
from $[0,1] \rightarrow {\mathbb R}$, endowed 
with the norm
$$
|f|_\beta=|f(0)|+ w_\beta (f, 1) \, 
$$
where
$$
w_\beta (f, x) =
\sup_{s,t \in [0,1], |t-s|\leq x} 
\frac{|f(t)-f(s)|}{|t-s|^\beta} \, .
$$
We denote by ${\mathcal H}^{0}_\beta([0,1])$
the subspace of ${\mathcal H}_\beta([0,1])$ 
of all functions $f$ such that 
$\lim_{x \rightarrow 0} w_\beta (f, x) =0$. It 
is well known (see \cite{C}) that ${\mathcal H}^{0}_\beta([0,1])$ is a closed subspace of   
${\mathcal H}_\beta([0,1])$, so that it
is a Banach space (for the norm $|\cdot|_\beta$),
and that it is separable (whereas 
${\mathcal H}_\beta([0,1])$ is not). 

Since the sample paths of the Brownian motion
are $\beta$-H\"older continuous for 
any $\beta < 1/2$, we may consider $W$ as 
a random variable taking values in 
${\mathcal H}^{0}_\beta([0,1])$, $\beta <1/2$. 
It is therefore natural  to look for sufficient conditions ensuring that  the convergence of $
\left \{ W_n(t), t \in [0,1]\right \}
$ to $\sigma W$ takes place in the space 
${\mathcal H}^{0}_\beta([0,1])$.  

In the case of strong mixing sequences in the sense of Rosenblatt, Giraudo  \cite{Gir} recently proved such an invariance principle under a sharp condition expressed  in terms of moments of the random variables and strong mixing rates. As we shall see, Giraudo's result can be extended to  $\alpha$-dependent sequences.

\begin{Theorem}\label{Holder}
Let $Q$ be an integrable quantile function, and let
  $Y_i$, $f$, $X_i$ and $S_n$ be as in Theorem \ref{ineRosalphafaible}. 
Assume that, for $p>2$,
\begin{equation}\label{WM0}
 \lim_{x \rightarrow \infty} x^{p-1} \int_0^1 Q(u) {\bf 1}_{R(u)>x} du  =0 \,  ,
\end{equation}
where $R$ is defined in \eqref{defofRuinfinity}. Let $\delta=(1/2)-(1/p)$. Then
the process
$
\left \{ W_n(t), t \in [0,1]\right \}
$
converges in distribution in ${\mathcal H}^0_{\delta}([0,1])$ to $\sigma W$, where $W$ is a standard 
Brownian motion and 
$
\sigma^2
$
is defined in \eqref{sigma}.
\end{Theorem}

\begin{remark}
In his paper \cite{Gir}, Giraudo obtains also sharp 
results for
$\tau$-dependent and $\rho$-mixing sequences. 
In a second paper \cite{Gir2}, he also proves
 H\"olderian
invariance principles
for other classes of stationary sequences $via$ martingale approximations.
\end{remark}

\begin{remark}
Applying Markov's inequality at order $p-1$, we see that the condition \eqref{SM} implies
\eqref{WM0}.
In the stationary $m$-dependent case,
we infer from Theorem \ref{Holder}  that the process
$
\left \{ W_n(t), t \in [0,1]\right \}
$
converges in distribution in ${\mathcal H}^0_{\delta}([0,1])$ to $\sigma W$
as soon as 
\begin{equation} \label{mdep2}
\lim_{x\rightarrow \infty}x^{p-1} {\mathbb E}(|X_0| {\bf 1}_{|X_0|>x}) =0 \, .
\end{equation} 
 Since $p>1$, It is easy to see (see Remark 2.3 in \cite{Gir})
that \eqref{mdep2} is equivalent to 
\begin{equation}\label{RS}
\lim_{x \rightarrow \infty}x^{p} {\mathbb P}(|X_0|>x) =0 \, .
\end{equation}
Note that, in the i.i.d. case, the condition
\eqref{RS} is necessary and sufficient for the invariance 
principle in ${\mathcal H}^0_{\delta}([0,1])$ (see
 \cite{RS}). 
\end{remark}

\section{Application to intermittent maps}

\subsection{Intermittent maps} \label{subsec:int}
Let us first recall the definition of the generalized Pomeau-Manneville maps introduced in 
\cite{DGM}. 

\begin{definition}
A map $\theta:[0,1] \to [0,1]$ is a generalized Pomeau-Manneville
map (or GPM map) of parameter $\gamma \in (0,1)$ if there exist
$0=y_0<y_1<\dots<y_d=1$ such that, writing $I_k=(y_k,y_{k+1})$,
\begin{enumerate}
\item The restriction of $\theta$ to $I_k$ admits a $C^1$ extension
$\theta_{(k)}$ to $\overline{I_k}$.
\item For $k\geq 1$, $\theta_{(k)}$ is $C^2$ on $\overline{I_k}$, and $|\theta_{(k)}'|>1$.
\item $\theta_{(0)}$ is $C^2$ on $(0, y_1]$, with $\theta_{(0)}'(x)>1$ for $x\in
(0,y_1]$, $\theta_{(0)}'(0)=1$ and $\theta_{(0)}''(x) \sim c
x^{\gamma-1}$ when $x\to 0$, for some $c>0$.
\item $\theta$ is topologically transitive.
\end{enumerate}
\end{definition}
The third condition ensures that $0$ is a neutral fixed point
of $\theta$, with $\theta(x)=x+c' x^{1+\gamma} (1+o(1))$ when $x\to 0$.
The fourth condition is necessary to avoid situations where
there are several absolutely continuous invariant measures, or
where the neutral fixed point does not belong to the support of
the absolutely continuous invariant measure.

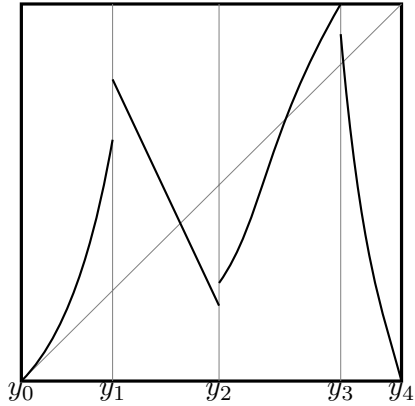
\begin{figure}[htb]
\centering
  \begin{tikzpicture}[scale=1]
  \draw[very thick] (0,0) rectangle (5,5);
  \draw[gray, very thin]
      (1.2,0) -- +(0, 5)
      (2.6,0) -- +(0, 5)
      (4.2,0) -- +(0, 5)
      (0,0)   -- (5,5);
  \draw[thick]
      (0,0) .. controls +(45:1) and +(-100:1) .. (1.2, 3.2)
      (1.2, 4) -- (2.6, 1)
      (2.6, 1.3) .. controls +(55:1) and +(-120:2) .. (4.2, 5)
      (4.2,4.6) .. controls +(-85:3) and +(105:1) .. (5, 0);
  \foreach \x/\ytext in {0/$y_0$, 1.2/$y_1$, 2.6/$y_2$, 4.2/$y_3$, 5/$y_4$}
      \node[above] at (\x, -0.4) {\ytext};
  \end{tikzpicture}
\caption{The graph of a GPM map, with $d=4$}
\end{figure}
The following  well known example of GPM map with only two branches 
has been introduced 
by  Liverani, Saussol and Vaienti \cite{LSV}:
\beq \label{LSVmap}
   \theta(x)=
  \begin{cases}
  x(1+ 2^\gamma x^\gamma) \quad  \text{ if $x \in [0, 1/2[$}\\
  2x-1 \quad \quad \quad \ \  \text{if $x \in [1/2, 1]$.}
  \end{cases}
\eeq

As quoted in \cite{DGM}, a GPM map $\theta$ admits 
a unique invariant absolutely continuous (with respect to the Lebesgue measure) probability  $\nu$ with density $h$. 
Moreover, it is ergodic, has full support, and $x^\gamma h(x)$ is bounded  above and below by positive constants.

We shall apply the  results of Sections \ref{Sec:MR}, \ref{Sec:LD} and \ref{Sec:HIP}  to the partial
sums
\begin{equation}\label{Bsums}
S_n(f)= \sum_{k=1}( f \circ \theta^k - \nu(f))
\end{equation}
where 
$\theta$ is a GPM map, and $f$ belongs to the space $\tMonm(Q, \nu)$
for 
some integrable quantile function $Q$. 

To do this, we shall make use of the results of \cite{DGM}. It is proved in that paper that there exists a stationary 
Markov chain ${\bf Y}=(Y_i)_{i\in {\mathbb Z}}$ such
that, on the probability space $([0, 1], \nu)$, the
random vector $(\theta, \theta^2, \ldots , \theta^n)$ is distributed as
$(Y_n,Y_{n-1}, \ldots, Y_1)$. 

 In particular, on $([0, 1], \nu)$, the partial sum
$S_n(f)$ defined in \eqref{Bsums} is distributed as $\sum_{k=1}^n X_i$ with
$X_i=f(Y_i)-\nu(f)$.
Moreover, since 
\begin{equation*}
  \label{equ1law}
  \max_{1 \leq k \leq
  n}  |S_k(f)|  \quad \text{is distributed as} \quad  \max_{1 \leq k \leq
  n} \left| \sum_{i=k}^n X_i \right | \, ,
\end{equation*}  
we easily derive that,
 for any $x\geq 0$,
\begin{equation}
\label{equ2law}
\nu \left( \max_{1 \leq k \leq n} |S_k(f)| \geq x \right) \leq
{\mathbb P}\left(2\max_{1 \leq k \leq n} \left | \sum_{k=1}^n X_i \right | \geq x \right) \,  .
\end{equation}



From Proposition 1.17 (and the comments right after)  in \cite{DGM}, we know that
for any positive integer $k$, there exist two  positive constants $C$ and $D$ such
that, for any $n>0$,
\[
 \frac{D}{n^{(1-\gamma)/\gamma}} \leq \alpha_{k,{\bf Y}}(n) \leq
  \frac{C}{n^{(1-\gamma)/\gamma}} \, .
\]
This control of the coefficients $\alpha_{k,{\bf Y}}(n)$
(for $k=2$), together with Inequality \eqref{equ2law} and the control 
 $\nu(|f|>t)\leq H(t)$, are all  we need to apply the results of Sections \ref{Sec:MR}, \ref{Sec:LD}
 and \ref{Sec:HIP}
to the partial sums $S_n(f)$.

\subsection{Moment bounds for intermittent maps}
\label{MBI}
In this subsection, we shall prove moment inequalities 
for $\max_{1 \leq k \leq n} |S_k(f)|$ when 
$f \in \tMonm(Q, \nu)$ and $Q(u) \ll u^{-b}$ for 
$b \in [0,1)$.

In particular, our results apply to the  two simple examples:
\begin{enumerate}
\item If $f$ is positive and non increasing on (0, 1), with
\begin{equation}\label{blow0}
f(x) \leq \frac{C}{x^{s}}
\quad \text{near 0, for some $C>0$ and  $s\in [0,  1-\gamma)$,}
\end{equation}
then $f$ belongs to $\tMonm(Q, \nu)$ with $Q(u)\ll u^{-s/(1-\gamma)}$  (meaning that 
$b=s/(1-\gamma)$). 
\item
If $f$ is positive and non increasing on (0, 1), with
\begin{equation}
f(x) \leq \frac{C}{(1-x)^{s}}
\quad \text{near 1, for some $C>0$ and  $s\in [0,  1)$,}
\end{equation}
then $f$ belongs to $\tMonm(Q, \nu)$ with $Q(u)\ll u^{-s}$  (meaning that 
$b=s$). 
\end{enumerate}

$\bullet$ Assume first that $ p >2$ and 
$b \in [0, 1/p)$ in such a way that 
$\int_0^1 Q^p (u) du < \infty$ for  $p>2$. 
From \eqref{equ2law} and 
Theorem \ref{ineRosalphafaible}, we infer that, for any $f$ in $\tMonm(Q, \nu)$,
\begin{equation}
\left \| \max_{1 \leq k \leq n} |S_k(f)|  \right \|_p^p
\ll
\begin{cases}
   n^{p/2} \quad \quad \quad  \ \ \quad \quad \quad \  \text{if $b\leq (2-\gamma(p+2))/(2p(1-\gamma))$}\\
   n^{(p \gamma+(\gamma-1)(1-pb))/\gamma}  \quad  
  \text{if $ b >(2-\gamma(p+2))/(2p(1-\gamma)) $.}
  \end{cases}
\end{equation}

$\bullet$ Assume now that $p=2$ and $b \in [0, 1/2)$ 
in such a way that 
$\int_0^1 Q^2 (u) du < \infty$. From \eqref{equ2law} and 
Theorem \ref{ineRosalphafaible}, we infer that, for any $f$ in $\tMonm(Q, \nu)$,
\begin{equation}\label{p=2}
\left \| \max_{1 \leq k \leq n} |S_k(f)|  \right \|_2^2
\ll
\begin{cases}
   n \quad \quad \quad  \ \, 
   \quad \quad  \quad  \quad \quad \text{if $b< (1-2\gamma)/(2(1-\gamma))$}\\
   n \ln (n) \quad  \quad \quad \, 
   \quad \quad \quad \text{if $b= (1-2\gamma)/(2(1-\gamma))$} \\
   n^{(2 \gamma+(\gamma-1)(1-2b))/\gamma}  \quad  
\   \text{if $ b >(1-2\gamma)/(2(1-\gamma)) $.}
  \end{cases}
\end{equation}

$\bullet$ Assume now that $p \in (1,2)$ and 
$b\in [0, 1/p)$  in such a way that 
$\int_0^1 Q^p (u) du < \infty$ for  $p<2$. 

 Applying   Remark 8 in \cite{DM} (see section 4.2 in  \cite{DGM} for its application to $\alpha$-dependent sequences)
the following upper bounds hold.

\beq 
\left \| \max_{1 \leq k \leq n} |S_k(f)|  \right \|_p^p
\ll
\begin{cases}
   n \quad \quad \quad   \quad  \quad \quad
   \quad \quad  \ \ \,  \text{if $b<(1-p\gamma)/(p(1-\gamma))$}\\
 n  \ln (n)  \ \quad \quad \quad \quad \quad  
 \quad \,
 \text{if $b=(1-p\gamma)/(p(1-\gamma))$} \\
   n^{(p \gamma+(\gamma-1)(1-pb))/\gamma}  \quad  \ \
  \text{if $ b > (1-p\gamma)/(p(1-\gamma))$.}
  \end{cases}
  \nonumber
\eeq 

Moreover, if $b=(1-p\gamma)/(p(1-\gamma))$, Proposition 1 in \cite{DM}  gives the upper bound
\begin{equation}\label{better}
\nu \left (  \max_{1 \leq k \leq n} |S_k(f)|  \geq   x \right )
\ll \frac{n}{ x^p} \, .
\end{equation}

\begin{remark} 
When $f$ is a bounded variation function (which corresponds to the case $b=0$) all the rates above are given in Theorem 4.5 of \cite{DM2}, and these rates are optimal (see the discussion in Section 4.4 of \cite{DM2}). 
\end{remark}

\begin{remark}
Let us consider the case where $f$ satisfies \eqref{blow0}
 (in which case $b=s/(1-\gamma)$).

 In that case, if $s=(1-2 \gamma)/2$, the upper bound 
 \eqref{p=2}
  gives that 
$$
 \left \| \max_{1 \leq k \leq n} |S_k(f)|  \right \|_2^2
\ll 
n \ln (n) \, .
$$
This upper bound is in accordance with a result by Gou\"ezel 
\cite{Gou}.  He  proved that, if 
$f$ is exactly of the form
$f(x)=x^{-(1-2\gamma)/2}$ and $\theta$ is the LSV map defined by \eqref{LSVmap}, then $S_n(f)/{\sqrt{n \ln(n)}}$ converges in 
distribution to a non-degenerate Gaussian random variable.

Now, if $s=(1-p \gamma)/p$,  the upper bound \eqref{better}  holds. 
This  is also in accordance with a result by Gou\"ezel \cite{Gou}.  He  proved that, if 
$f$ is exactly of the form
$f(x)=x^{-(1-p\gamma)/p}$ and $\theta$ is the LSV map defined by \eqref{LSVmap}, then for any 
positive real $x$,
$$
\lim_{n \rightarrow \infty} \nu
\left(\frac{1}{n^{1/p}}\left | S_n(f)\right| > x \right)=
{\mathbb P}(|Z_p|>x) \, ,
$$
where $Z_p$ is a $p$-stable random variable such that $\lim_{x \rightarrow \infty}
x^p{\mathbb P}(|Z_p|>x)=c>0$.
\end{remark}

\subsection{Large deviations for intermittent maps}
Let $f$ in $\tMonm(Q, \nu)$ with $Q(u) \ll u^{-b}$ for some $b \in [0,1)$. 

$\bullet$ Assume first that 
$\gamma + b (1-\gamma) < 1/2$, and let 
$p=1/(\gamma + b (1-\gamma))$. 
Applying Inequality \eqref{WB}, we get that, 
for any  $a \in (p-1, p)$ and any $x>0$, 
\begin{equation*}
\nu \left(\frac 1 n  \max_{ 1 \leq k \leq n}  |
  S_k(f) | \geq x \right)
 \ll \frac 1 {n^a x^{2a}} + \frac 1 {n^{p-1} x^p} \, .
\end{equation*}

$\bullet$ Assume now that 
$\gamma + b (1-\gamma) = 1/2$. Applying Inequality 
\eqref{WB2} 
we get that,  
for any  $a \in (1, 2)$, any $c \in (0,1)$ and any $x>0$, 
\begin{equation*}
\nu \left(\frac 1 n  \max_{ 1 \leq k \leq n}  |
  S_k(f) | \geq x \right)
 \ll \frac {1} {n^{a c } x^{a(1+c )}} + \frac 1 {n x^2} \, .
\end{equation*}

$\bullet$ Assume now that 
$\gamma + b (1-\gamma)  \in (1/2, 1)$, and let 
$p=1/(\gamma + b (1-\gamma))$. 
 Applying  \eqref{better}, 
we get that,  
for  any $x>0$, 
\begin{equation*}
\nu \left(\frac 1 n  \max_{ 1 \leq k \leq n}  |
  S_k(f) | \geq x \right)
 \ll  \frac 1 {n^{p-1} x^p} \, .
\end{equation*}

\begin{remark}
Let $b \in [0, 1)$ and $p=1/(\gamma + b (1-\gamma))$.
From the preceeding upper bounds, we infer that there exists a function $f_{b, \gamma}$ from ${\mathbb R}^+$ to 
${\mathbb R}^+$
 such that
for any $x>0$,
\begin{equation}\label{summary}
\nu \left(\frac 1 n  \max_{ 1 \leq k \leq n}  |
  S_k(f) | \geq x \right)
 \ll  \frac {f_{b, \gamma}(x)} {n^{p-1}} \, .
\end{equation}
Moreover $\sup_{x >\varepsilon} x^p f_{b, \gamma}(x) < \infty$ for 
any $\varepsilon>0$.  
\end{remark}
\begin{remark} 
When $f$ is a bounded variation function (which corresponds to the case $b=0$), we  obtain from 
\eqref{summary} that, 
for any $x>0$,
\begin{equation}\label{Mel}
\nu \left(\frac 1 n  \max_{ 1 \leq k \leq n}  |
  S_k(f) | \geq x \right)
 \ll  \frac {f_{0, \gamma}(x)} {n^{(1-\gamma)/\gamma}} \, .
\end{equation}
Note that the upper bound \eqref{Mel} (with $S_n(f)$ instead of the maximum) has been already obtained by 
Melbourne (\cite{M}, Example 1.6)
when $\theta$  is the LSV map defined by \eqref{LSVmap} 
and $f$ is  H\"older continuous (as a consequence of a very  general result on slowly mixing dynamical system). In that case, the 
bound is optimal (see again Example 1.6 in  \cite{M}). 
\end{remark}

\subsection{H\"olderian invariance principles for intermittent maps}
Assume now 
that 
$f$ belongs to  $\tMonm(Q, \nu)$ for some $Q$ such that $Q(u) \leq u^{-b} \varepsilon(u)$, where
$b \in (0,1)$ and $\varepsilon$ is a bounded function such that $\varepsilon(u) \rightarrow 0$  as 
$u \rightarrow 0$.

In particular, our results apply to the  two simple examples:
\begin{enumerate}
\item If $f$ is positive and non increasing on (0, 1), with
\begin{equation}\label{blow0bis}
f(x) \leq \frac{\varepsilon(x)}{x^{s}}
\quad \text{near 0, for some   $s\in [0,  1-\gamma)$,}
\end{equation}
then the assumption on $Q$ is satisfied with
$b=s/(1-\gamma)$. 
\item
If $f$ is positive and non increasing on (0, 1), with
\begin{equation}
f(x) \leq \frac{\varepsilon(1-x)}{(1-x)^{s}}
\quad \text{near 1, for some   $s\in [0,  1)$,}
\end{equation}
then 
the assumption on $Q$ is satisfied with
$b=s$. 
\end{enumerate} 

Let then
$$
W_n(f,t)= \frac 1 {\sqrt{n}} S_{[nt]}(f) +\frac{(nt-[nt])}{\sqrt n} 
f \circ \theta^{[nt]+1} \, .
$$

Assume that
$\gamma + b (1-\gamma) < 1/2$,
and let $\delta= (1/2)- \gamma -b(1-\gamma)$.
Applying 
Theorem \ref{Holder}, we infer that, on the probabilty space $([0,1], \nu)$,  the process
$
\left \{ W_n(f, t), t \in [0,1]\right \}
$
converges in distribution in ${\mathcal H}^0_{\delta}([0,1])$ to $\sigma(f) W$, where $W$ is a standard 
Brownian motion and 
$$
\sigma^2(f)= \nu\left((f-\nu(f))^2\right)
+
2 \sum_{k>0} \nu\left((f-\nu(f))\cdot  f \circ T^k\right)
\, .
$$

\begin{remark} 
When $f$ is a bounded variation function and $\gamma < 1/2$, we infer 
that  the process $
\left \{ W_n(f, t), t \in [0,1]\right \}
$ converges in distribution in ${\mathcal H}^0_{\delta}([0,1])$ to $\sigma(f) W$, for any $\delta<1/2-\gamma$. 
\end{remark}

\section{A  deviation inequality for the maximum of partial sums} \label{sectiondeviation}

Recall first that the functions $\alpha_{1,{\bf Y}}^{-1}$ and  $\alpha_{2,{\bf Y}}^{-1}$ have been defined in \eqref{invalpha}. Our key deviation inequality for $\alpha$-dependent sequences is given below. 

\begin{Proposition} \label{deviationRosenthal}
 Let ${\bf Y}=(Y_i)_{i \in {\mathbb Z}}$,
where $Y_i= Y_0 \circ T^i$.
Denote by $P_{Y_0}$  the distribution of $Y_0$, and let $Q$ be a
 quantile function in ${\mathbb L}^1$.
Let  $X_i = f(Y_i) - \E ( f(Y_i))$, where  $f$ belongs to $\tMonm(Q, P_{Y_0})$.  Given a positive integer $n$, define
  \begin{equation*}
  R_n(u)=\left(  \alpha_{2,{\bf Y}}^{-1} (u)\wedge n\right ) Q(u)\, ,\quad 
  \text{for $ u  \in [0,1]$}
  \end{equation*}
   and 
 \begin{equation*}
  L_n(x)=R_n^{-1}(x)=\inf
  \left \{u\in [0,1] \tq R_n(u)\leq x \right \} \, ,\quad 
  \text{for $x \geq 0$} \, .
  \end{equation*}
Let $S_n = \sum_{k=1}^n X_k$. For any $x>0$, $r >2$, $\beta \in ]r-2,r[$ and 
$$
s_n^2(x) = n  \int_{L_n(x)}^1 (  \alpha_{1, {\bf Y}}^{-1} (u) \wedge n)
Q^2(u)du \, ,
$$ 
the following deviation bound holds
\begin{multline}
  \label{eqmainmax}
  \p \left( \max_{ 1 \leq k \leq n}  |
  S_k | \geq x \right)  \ll 
  \frac{s^r_n(x)}{x^r}+   \frac{n}{x} \int_0^{L_n(x)}Q(u) du +  \frac{n}{x^{1+\beta/2}}   \int_0^{L_n(x)} R_n^{\beta/2}(u)   Q(u) du \\
  +   \frac{n}{x^{1+r/2}}   \int_{L_n(x)}^1 R_n^{r/2}(u)   Q(u) du \, .
 \end{multline}
\end{Proposition}
\begin{remark} \label{rmk:sn2}
The quantity $s_n^2(x)$ can be bounded as follows:
\begin{equation}\label{sn2}
s_n^2(x) \leq s_n^2 :=n \sum_{k=0}^{n-1} \int_0^{\alpha_{1,
{\bf Y}} (k)}Q^2 (u) du \, .
\end{equation}
\end{remark}
\noindent {\bf Proof of Proposition \ref{deviationRosenthal}}
In all the proof, we shall use the following notation:
for any non-negative random variable $V$, let $Q_V$ be the generalized inverse of the tail function $x \rightarrow \p (V>x)$. 

We shall first prove Proposition \ref{deviationRosenthal} 
for $X_i = \sum_{\ell=1}^L a_{\ell} f_{\ell}(Y_i)
- \sum_{\ell=1}^L a_{\ell}\E(f_{\ell}(Y_i))$, with $f_\ell$
belonging to $\tMon(Q, P_{Y_0})$ and $\sum_{\ell=1}^L |a_\ell|
\leq 1$. 
Let 
\beq  \label{selectvandM}
v = L_n(x) \,  \text{ and }  \, M=Q(v) \, .
\eeq 
Since $\alpha_{2, {\bf Y}} (0) =1/2$, it follows that  $v \in [0,1/2]$.  Note first that if $v =1/2$, then by using Markov's inequality, we derive
\begin{multline*}
\p \big ( \max_{1 \leq k \leq  n } |S_k |  \geq x \big )  \leq  \frac{1}{x} \sum_{k=1}^n \E (| X_k | ) \leq  \frac{2}{x}   \sum_{\ell=1}^L | a_{\ell} |  \sum_{k=1}^n \E ( |  f_{\ell}(Y_k) |)  = \frac{2 n }{x}   \sum_{\ell=1}^L | a_{\ell} |  \int_0^1 Q_{| f_{\ell}(Y_k)|}(u) du   \\
\leq \frac{2 n }{x}    \int_0^1 Q(u) du   \leq \frac{4 n }{x}    \int_0^{1/2} Q(u) du = \frac{4 n }{x}    \int_0^{L_n(x)} Q(u) du    \, , 
\end{multline*}
which then proves the proposition in case where $v=1/2$. 

Therefore, we can assume in the rest of the proof that $v < 1/2$.  We then set $ g_M(y) = (y \wedge M) \vee (-M) $ where $M$ is defined in \eqref{selectvandM}, and, for
any $i \in {\mathbb Z}$ and any $\ell \geq 1$, we  define
\[
Y_{i,\ell}'= g_M \circ f_{\ell}(Y_i) - \E ( g_M \circ f_{\ell}(Y_i)) \quad \text{and} \quad  Y_{i,\ell}''=f_{\ell}(Y_i)  - 
\E(f_{\ell}(Y_i)) - Y_{i,\ell}'  \, .
\]
Then, for any $i \in {\mathbb Z}$, we set
\[
 X_i'=\sum_{\ell=1}^L a_{\ell} Y_{i,\ell}' \quad \text{and} \quad  X_i''=X_i - X_i' \, .
\]
Let $q = \min\{k\in {\mathbb N} \tq \alpha_{2,{\bf
Y}}(k)\leq v\} \wedge n $ where $v$ in defined in \eqref{selectvandM}.  Since $R_n$ is right-continuous and non-increasing, we have $R_n(L_n(w))\leq w$ for any $w$, hence
\beq \label{restq}
qM = R_n(v) = R_n (L_n(x)) \leq x  \, .
\eeq
Assume first that $q=n$. Bounding  $X'_i$ by $2M$, we obtain
$\max_{1\leq k\leq n}|S_k| \leq 2qM + \sum_{k=1}^n |X''_k|$.
Taking into account \eqref{restq} this gives
  \begin{equation*}
  \p \left ( \max_{1 \leq k \leq  n } |S_k | \geq 8x \right)
  \leq \frac{1}{6x} \sum_{k=1}^n  \E ( | X_k''|).
  \end{equation*}
Writing $\varphi_M(x)=(|x|-M)_+$, we have
\[
\sum_{k=1}^n \E ( | X_k''|) \leq 2 \sum_{\ell =1}^L
|a_{\ell}|\sum_{k=1}^n \E (\varphi_{M} ( f_{\ell} (Y_k))) \, .
\]
But $Q_{\varphi_M ( f_{\ell} (Y_k))} \leq Q_{| f_{\ell}
(Y_k) |}\I_{[0,v]} \leq Q\I_{[0,v]} $. Consequently
  \begin{equation} \label{dec13FN}
  \sum_{k=1}^n  \E ( | X_k''|)
  \leq  2 \sum_{\ell =1}^L |a_{\ell}|\sum_{k=1}^n  \int_0^v  Q_{|
  f_{\ell} (Y_k) |} (u) du
  \leq  2n \int_0^v Q(u) du \leq 2n \int_0^{L_n(x)} Q(u) du\, ,
  \end{equation}
  where for the last inequality we have used that  $u<v \iff Q(v) <Q(u)$. We derive that
\begin{equation*}
   \p \left ( \max_{1 \leq k \leq  n } |S_k | \geq 8x \right)
  \leq  \frac{n}{x} \int_0^{L_n(x)} Q(u) du \, ,
  \end{equation*}
which proves the result in case where $q=n$.

From now on, we assume that $q<n$. Therefore $q =  \min\{k\in {\mathbb N} \tq \alpha_{2,{\bf
Y}}(k)\leq v\} $ and then $\alpha_{2,{\bf
Y}}(q) \leq v$.  Recall also that since $v$ is assumed to be strictly less than $1/2$ then  $q \geq 1$.
For any  integer $i$, define the random variables
\[ U_i =  \sum_{k=(i-1)q+1}^{iq} X_k'  \, .\]
Consider now the $\sigma$-algebras ${\cal G}_{i}= \F_{iq}$ and define the  variables $\tilde
U_i$ as follows:   
$$\tilde U_{2i -1} =  U_{2i -1}-
\E (U_{2i - 1 } | {\cal G}_{2(i-1) -1 }) \quad 
\text{and} \quad  \tilde U_{2i}=U_{2i}- \E (U_{2i} | {\cal
G}_{2(i-1)}) \, .
$$ 
The following inequality is then valid
\begin{multline*}
\max_{1 \leq k \leq  n } |S_k |\leq  2q M +
 \max_{2 \leq 2j \leq [n/q]} \left | \sum_{i=1}^j \tilde U_{2i} \right | +
\max_{1 \leq 2j-1 \leq [n/q]} \left | \sum_{i=1}^j \tilde U_{2i
-1} \right |
  + \sum_{i=1}^{ [n/q]} |  U_{i} -\tilde U_{i} | +  \sum_{k=1}^n
| X''_k |  \, .
\end{multline*}
(See the proof of Inequality (2.2) in \cite{DGM}). Using \eqref{dec13FN} and \eqref{restq}, it follows that 
\begin{multline} \label{dec15FN}
\p \left ( \max_{1 \leq k \leq  n } |S_k | \geq 8x \right ) \leq  
 \p \left (\max_{2 \leq 2j \leq [n/q]} \left | \sum_{i=1}^j \tilde U_{2i} \right | \geq x \right )+
\p \left (\max_{1 \leq 2j-1 \leq [n/q]} \left | \sum_{i=1}^j \tilde U_{2i
-1} \right | \geq  x \right ) \\
  + \p \left ( \sum_{i=1}^{ [n/q]} |  U_{i} -\tilde U_{i} | \geq x 
  \right ) +  \frac{n}{x} \int_0^{L_n(x)} Q(u) du \, .
\end{multline}
Using Markov's inequality and stationarity, we get 
\[
 \p \left ( \sum_{i=1}^{ [n/q]} |  U_{i} -\tilde U_{i} | \geq  x  \right
 ) \\
    \leq \frac{1}{x} \sum_{i=1}^{[n/q]} \left \| \E (  U_{i}  | {\cal
F}_{(i-2)q} ) \right \|_1 = \frac{n}{q x} \sum_{i=q+1}^{2q}
\left \| \E ( X'_i | {\cal F}_{0} ) \right \|_1\, .
\]
Now, setting $A = \sign \{ \E ( X'_i | {\cal F}_{0} ) \}$, 
\[
\| \E ( X'_i | {\cal F}_{0} )\|_1 = \E ((A - \E A)X_i')=
\sum_{\ell=1}^L a_{\ell} \E \left ( (A - \E A) (g_M \circ
f_{\ell}(Y_i) - \E g_M \circ f_{\ell}(Y_i))\right ) \, .
\]
Applying  Theorem 1.1 in  \cite{Ri},
 and using the fact that $Q_{|g_M \circ f_{\ell}(Y_i)|}
(u) \leq Q(u)$, we derive that
\begin{equation*}
 \left \vert  \E \left ( (A - \E A) (g_M \circ f_{\ell}(Y_i) - \E g_M
\circ f_{\ell}(Y_i))\right ) \right \vert
  \leq 2 \int_0^{ 2\bar \alpha ( A, g_M \circ
f_{\ell}(Y_i))} Q(u) du
 \, ,
\end{equation*}
where, for any real-valued random variables $A, B$,
$$
\bar \alpha ( A, B)=
\sup_{(s,t) \in {\mathbb R}^2} \left | {\mathrm{Cov}}
\left (
{\mathbf 1}_{A \leq s}, {\mathbf 1}_{B\leq t}
\right )
\right | \, .
$$
Since for all $i \geq q$,
\[
\bar \alpha ( A, g_M \circ f_{\ell}(Y_i) ) \leq 2 \bar
\alpha(A, Y_i) \leq  \alpha_{1, {\bf Y}} (i) \leq  \alpha_{2,
{\bf Y}} (i) \, ,
\]
we derive 
  \begin{equation}
  \label{majnorm1prime}
  \left \| \E ( X'_i | {\cal F}_{0} )
  \right \|_1 \leq 4 \int_0^{\alpha_{2, {\bf Y}} (i)} Q(u) du \, ,
  \end{equation}
which implies 
  \begin{equation}
  \label{decinter2}
  \p\left ( \sum_{i=1}^{ [n/q]} |  U'_{i} -\tilde U_{i} |\geq  x  \right
  ) \leq \frac{ 4 n}{x}\int_0^{\alpha_{2, {\bf Y}} (q)} Q (u) du \leq  \frac{ 4 n}{x}\int_0^{L_n(x)} Q (u) du
 \, .
  \end{equation}
Next, by using Markov's inequality, 
\beq \label{MarkovBigbloc}
 \p \left(\max_{2 \leq 2j \leq [n/q]} \left | \sum_{i=1}^j \tilde U_{2i} \right | \geq x \right ) \leq \frac 1 {x^{r}}
 \left \Vert\max_{2 \leq 2j \leq [n/q]} \left | \sum_{i=1}^j \tilde U_{2i} \right | \right \Vert_r^r\, .
\eeq
Note that $( \tilde U_{2i})_{i \in {\mathbb Z}} $ (resp. $( \tilde U_{2i -1
})_{i \in {\mathbb Z}}$) is a stationary sequence of martingale differences  with respect to the filtration $ ( {\cal {G}}_{2i} )_{i \in {\mathbb Z}}$ (resp. $ ( {\cal {G}}_{2i -1 } )_{i \in {\mathbb Z}}$). By using the Rosenthal inequality given in \cite{MP} (see their Theorem 6), we get
\[
\left\Vert \max_{2 \leq 2j \leq
 [n/q]
 } \left\vert \sum_{i=1}^j \tilde U_{2i}  \right\vert \right \Vert_r 
\ll (n/q)^{1/r}  \Vert \tilde U_{2}  \Vert_r + (n/q)^{1/r} 
\left (  \sum_{k=1}^{[n/q]} \frac{1}{k^{1+2\delta/r}}  
\left \Vert  \E_0 \left ( \left ( \sum_{i=1}^k \tilde U_{2i}  
\right )^2\right )\right \Vert_{r/2}^{\delta}\right )^{\frac 1{2\delta}}\, ,
\]
where $\delta= \min ( 1, 1/(r-2))$. Since $( \tilde U_{2i})_{i \in {\mathbb Z}} $ is a stationary sequence of martingale differences  with respect to the filtration $ ( {\cal {G}}_{2i} )_{i \in {\mathbb Z}}$, 
\[
\E_0 \left ( \left ( \sum_{i=1}^k \tilde U_{2i}  \right )^2 \right ) = \sum_{i=1}^k \E_0 \left (  \tilde U^2_{2i}    \right ) \, .
\]
Moreover, 
\[
\E_0 \left (  \tilde U^2_{2i}  \right ) \leq \E_0 
\left ( U^2_{2i}   \right ) \, .
\]
Therefore
\[
\left \Vert  \E_0 \left ( \left ( \sum_{i=1}^k \tilde U_{2i} 
 \right )^2\right )\right \Vert_{r/2} \leq \sum_{i=1}^k  \Vert 
 \E_0 \left ( U^2_{2i}  \right ) - \E \left (   U^2_{2i}  \right )\Vert_{r/2} +  \sum_{i=1}^k  \E \left (U^2_{2i} \right ) \, .
\]
By stationarity
\[
\sum_{i=1}^k  \E \left ( U^2_{2i}  \right  ) = k \Vert S'_q  \Vert_2^2 \, , \mbox{ with} \ S_q' = \sum_{i=1}^q X_i' \, .
\]
Moreover
\[
\left \Vert \tilde U_{2}  \right \Vert_r \leq 
 2 \left \Vert S_q' \right \Vert_r \, .
\]
It follows that
\beq \label{consRosenthal}
\left\Vert \max_{2 \leq 2j \leq
 [n/q]
 } \left\vert \sum_{i=1}^j \tilde U_{2i} (t) \right\vert \right\Vert^r_r \\
\ll \frac{n}{q} 
\left \Vert  S'_q \right  \Vert^r_r + 
\left ( \frac{n}{q} \right )^{r/2} 
\left \Vert S'_q  \right \Vert^r_2 +  \frac{n}{q} \left (  \sum_{k=1}^{[n/q]} \frac{1}{k^{1+2\delta/r}} D_{k,q}^{\delta}\right )^{r/(2\delta)}\, ,
\eeq
where 
\[
D_{k,q}  = \sum_{i=1}^k  
\left \Vert \E_0 \left ( U^2_{2i} \right ) - \E \left (   U^2_{2i}  \right ) \right \Vert_{r/2}\, .
\]
Notice that 
\begin{multline*}
\left \Vert \E_0 \left ( U^2_{2i} \right ) - \E \left (   U^2_{2i}  \right )\right \Vert_{r/2} \leq \sum_{j,k =(2i-1)q+1}^{2iq}  
\left \Vert \E_0  ( X'_j X'_k ) - \E ( X'_j X'_k ) \right \Vert_{r/2} \\
 \leq \sum_{\ell , m = 1}^L\sum_{j,k =(2i-1)q+1}^{2iq} |a_{\ell} | | a_m|  
 \left \Vert \E_0  ( Y'_{j,\ell} Y'_{k,m} ) - \E ( Y'_{j,\ell} Y'_{k,m} ) \right \Vert_{r/2} \, .
\end{multline*}
Now setting
\[
Z: =\left | \E_0  ( Y'_{j,\ell} Y'_{k,m} ) - \E ( Y'_{j,\ell} Y'_{k,m} )\right |^{r/2-1} \sign \left \{ \E_0  ( Y'_{j,\ell} Y'_{k,m} ) - \E ( Y'_{j,\ell} Y'_{k,m} ) \right \} \, \] 
we get 
\begin{multline*}
 \left \Vert \E_0  ( Y'_{j,\ell} Y'_{k,m} ) - \E ( Y'_{j,\ell} Y'_{k,m} ) \right \Vert_{r/2}^{r/2}
  =  \bkE \left ( Z  \left (\E_0  ( Y'_{j,\ell} Y'_{k,m} ) - \E ( Y'_{j,\ell} Y'_{k,m} ) \right )\right ) \\
  =  \bkE \left ( (Z- \bkE(Z))  Y'_{j,\ell} Y'_{k,m} \right )  \, .
  \end{multline*}
From Proposition A.1 and Lemma A.1 in \cite{DR2}, noticing that $Q_{|g_M \circ f_{\ell}(Y_i)|} (u) \leq M$, we derive
\begin{eqnarray*}
\bkE \left ( (Z- \bkE(Z))  Y'_{j,\ell} Y'_{k,m} \right )  \leq 2^{4} M^2  \int_0^{\bar \alpha /2} Q_{|Z|}(u) du\, ,
\end{eqnarray*}
with $\bar \alpha:=  \alpha ( Z, g_M \circ f_{\ell}(Y_j), g_M \circ f_{m}(Y_k) )$,
where for real-valued random variables $Z,B,V$,
\[
  \bar \alpha (Z, B,V) = 
  \sup_{(s,t,u) \in {\mathbf R}^3} \left | \E ( (\I_{ Z
  \leq s} - \p (Z \leq s)) (\I_{ B \leq t} - \p (B \leq t)) (\I_{ V
  \leq u} - \p (V \leq u)))\right | \, .
\]
Applying H\"older's inequality,
\[
\int_0^{\bar \alpha /2} Q_{|Z|}(u) du \leq (\bar \alpha /2 )^{2/r} 
\left ( \int_0^1 Q_{|Z|}^{r/(r-2)} (u)\right )^{(r-2)/r} = (\bar \alpha /2 )^{2/r} \Vert Z \Vert_{r/(r-2)} \, .
\]
Since $$\Vert Z \Vert_{r/(r-2)} =  
\left \Vert \E_0  ( Y'_{j,\ell} Y'_{k,m} ) - \E ( Y'_{j,\ell} Y'_{k,m} ) \right \Vert_{r/2}^{r/2 -1} \, , $$
 we then derive
\[
\left \Vert \E_0  ( Y'_{j,\ell} Y'_{k,m} ) - 
\E ( Y'_{j,\ell} Y'_{k,m} ) \right \Vert_{r/2} 
\leq 2^{4} M^2 (\bar \alpha /2 )^{2/r} \, .
\]
Now, for all $j,k \geq (2i-1)q+1$,
\[
\bar \alpha  \leq 4 \bar \alpha (Z, f_{\ell}(Y_j),f_{\ell}(Y_k)) \leq 2 \alpha_{2, {\bf Y}}((2i-q) +1) \, .
\]
So, overall, using the fact that $\sum_{\ell = 1}^L |a_{\ell} | \leq 1$ and that $\alpha_{2, {\bf Y}}((2i-q) +1) \leq \alpha_{2, {\bf Y}}(iq+1)$, we get 
\[
D_{k,q}  \leq 2^{4} q^2 M^2 \sum_{i=1}^k \alpha_{2, {\bf Y}}^{2/r}(iq+1)\, .
\]
Let $\eta = (\beta-2)/r$ and recall that $
r-2 < \beta < r$. 
Since $\eta< (r-2)/r$, applying H\"older's inequality, we then get 
\[
D_{k,q} \ll q^2 M^2 k^{- \eta +(r-2)/r } \left ( \sum_{i=1}^k i^{\beta/2-1}\alpha_{2, {\bf Y}} (iq+1)  \right )^{2/r} \, .
\]
Since $2\delta/r > - \delta \eta +\delta (r-2)/r$ (indeed $ - \eta + (r-2)/r  = (r-\beta)/r$ and $r-\beta < 2$), it follows that 
\[
\sum_{k=1}^{[n/q]} \frac{1}{k^{1+2\delta/r}} D_{k,q}^{\delta} \ll (qM)^{2\delta}  \left ( \sum_{i=1}^{[n/q]} i^{\beta/2-1}\alpha_{2, {\bf Y}} (iq+1)  \right )^{2 \delta/r} \, .
\]
But, since $x < \alpha_{2, {\bf Y}}^{-1} (u) \iff \alpha_{2, {\bf Y}}(x) >u$ and $\alpha_{2, {\bf Y}}(q) \leq v$ (since $q<n$),
 \begin{multline*}
\sum_{i=1}^{[n/q]} i^{\beta/2-1}\alpha_{2, {\bf Y}} (iq+1) = \sum_{i=1}^{[n/q]} i^{\beta/2-1} \int_0^1 {\bf 1}_{u < \alpha_{2, {\bf Y}} (iq+1)} du \\
 \leq  \int_0^v \sum_{i=1}^{[n/q]} i^{\beta/2-1}   {\bf 1}_{i \leq q^{-1} \alpha_{2, {\bf Y}}^{-1} (u)}  du
 \leq   q^{-\beta/2}\int_0^v ( \alpha_{2, {\bf Y}}^{-1} (u) \wedge n)^{\beta/2} du \, .
\end{multline*}
So, overall,  
\[
\sum_{k=1}^{[n/q]} \frac{1}{k^{1+2\delta/r}} D_{k,q}^{\delta} \ll (qM)^{2\delta} q^{-\beta \delta /r} \left ( \int_0^v ( \alpha_{2, {\bf Y}}^{-1} (u) \wedge n)^{\beta/2} du \right )^{2 \delta/r} \, ,
\]
which implies that 
\[
\frac{n}{qx^r} \left (  \sum_{k=1}^{[n/q]} \frac{1}{k^{1+2\delta/r}} D_{k,q}^{\delta}\right )^{r/(2\delta)} \ll  \frac{n}{x^r} (qM)^{r} q^{-\beta /2 -1 }  \int_0^v ( \alpha_{2, {\bf Y}}^{-1} (u) \wedge n)^{\beta/2} du \, .
\]
Using the fact that $r-\beta/2-1 >0$ (since $\beta < r $ and $r>2$) and \eqref{restq}, we infer that 
$(qM)^{r} q^{-\beta /2 -1 } \leq x^{r-\beta /2 -1 } Q^{\beta/2 +1} (v)$. Moreover, $u<v \iff Q(v) <Q(u)$. Therefore
\beq \label{majorationDkq}
\frac{n}{qx^r} \left (  \sum_{k=1}^{[n/q]} \frac{1}{k^{1+2\delta/r}} D_{k,q}^{\delta}\right )^{r/(2\delta)} \ll  \frac{n}{x^{\beta/2+1}}  \int_0^v ( \alpha_{2, {\bf Y}}^{-1} (u) \wedge n)^{\beta/2} Q^{\beta/2+1} (u) du \, .
\eeq
We prove now that 
\beq \label{majorationSqr}
 \frac{n}{qx^r} \Vert  S'_q  \Vert^r_r \ll \frac{n}{x} \int_0^{L_n(x)}Q(u) du 
  +   \frac{n}{x^{1+r/2}}   \int_{L_n(x)}^1 R_n^{r/2}(u)   Q(u) du \,  .
\eeq
Using the stationarity and  applying Theorem 2.5(a) in \cite{Ri}, we get 
\beq \label{TH25RIO}
\Vert  S'_q  \Vert_r \leq \sum_{\ell =1}^L |a_{\ell} | \left \Vert \sum_{i=1}^q Y'_{i,\ell}\right \Vert_r \leq \sum_{\ell =1}^L |a_{\ell} | \left ( 2  q r   \max_{1 \leq k \leq q} \left \Vert Y_{0,\ell}' \sum_{i=0}^{k-1}\E_0 (Y_{i,\ell}') \right \Vert_{r/2}\right )^{1/2} \, .
\eeq
Setting
\[
W_{0,k}: =\left | Y_{0,\ell}'\sum_{i=0}^{k-1}\E_0 (Y_{i,\ell}') 
\right |^{r/2-1} 
\sign \left \{ Y_{0,\ell}' \sum_{j=0}^{k-1}\E_0  (Y'_{j,\ell}) \right \} \, \] we get
\[
 \left \Vert Y_{0,\ell}' \sum_{i=0}^{k-1}\E_0 (Y_{i,\ell}') \right \Vert^{r/2}_{r/2}= \sum_{j=0}^{k-1} \bkE \left ( W_{0,k} Y_{0,\ell}' \E_0  ( Y'_{j,\ell}) \right )   = \sum_{j=0}^{k-1} \bkE \left ( (W_{0,k}Y_{0,\ell}'- \bkE(W_{0,k} Y_{0,\ell}'))  g_M \circ f_{\ell}(Y_j) \right )  \, .
\]
From Proposition A.1 and Lemma A.1 in \cite{DR2}, noticing that $Q_{|g_M \circ f_{\ell}(Y_j)|} (u) \leq Q(u \vee v)$, we derive
\begin{multline*}
\bkE \left ( (W_{0,k}Y_{0,\ell}'- \bkE(W_{0,k}Y_{0,\ell}'))  g_M \circ f_{\ell}(Y_j) \right ) \leq 4  \int_0^{\tilde \alpha /2} Q_{|W_{0,k}Y_{0,\ell}'|}(u) Q(u \vee v) du \\ \leq 8  \int_0^{\tilde \alpha /2} Q_{|W_{0,k}|}(u) Q^2(u \vee v) du\, ,
\end{multline*}
where 
\[
  \tilde \alpha:= \sup_{(s,t) \in {\mathbf R}^2} \big | \E ( (\I_{ W_{0,k}Y_{0,\ell}'
  \leq s} - \p (W_{0,k}Y_{0,\ell}' \leq s)) (\I_{ g_M \circ f_{\ell}(Y_j) \leq t} - \p ( g_M \circ f_{\ell}(Y_j) \leq t)) )\big | \, .
\]
Now, for any $j\geq 0$,
\[
\bar \alpha \leq 2 \sup_{(s,t) \in {\mathbf R}^2} \big | \E ( (\I_{ W_{0,k}Y_{0,\ell}'
  \leq s} - \p (W_{0,k}Y_{0,\ell}'\leq s)) (\I_{f_{\ell}(Y_j) \leq t} - \p ( f_{\ell}(Y_j) \leq t)) ) \big | \leq 2  \alpha_{1, {\bf Y}}(j) \leq 2  \alpha_{2, {\bf Y}}(j) \, .
\]
It follows that 
\[
\left \Vert Y_{0,\ell}' \sum_{i=0}^{k-1}\E_0 (Y_{i,\ell}') \right \Vert^{r/2}_{r/2} \leq 8 \sum_{j=0}^{k-1} \int_0^1 \I_{u < \alpha_{2, {\bf Y}}(j)} Q_{|W_{0,k}|}(u) Q^2(u \vee v) du  \, .
\] 
Notice that $\sum_{j=0}^{k-1} \I_{u < \alpha_{2, {\bf Y}}(j)} = \alpha_{2, {\bf Y}}^{-1} (u) \wedge k$. Applying H\"older's inequality, we then get 
\begin{multline*}
\sum_{j=0}^{k-1} \int_0^1 \I_{u < \alpha_{2, {\bf Y}}(j)} Q_{|W_{0,k}|}(u) Q^2(u \vee v) du \\ \leq \left ( \int_0^1 (  \alpha_{2, {\bf Y}}^{-1} (u) \wedge k)^{r/2} Q^{r} (u \vee v) du \right )^{2/r}\left ( \int_0^1 Q_{|W_{0,k}|}^{r/(r-2)} (u)du \right )^{(r-2)/r} \\
\leq \left ( \int_0^1 (  \alpha_{2, {\bf Y}}^{-1} (u) \wedge k)^{r/2} Q^{r} (u \vee v) du \right )^{2/r}   \left \Vert  Y_{0,\ell}'   \sum_{i=0}^{k-1}\E_0 (Y_{i,\ell}') \right \Vert^{r/2-1}_{r/2} \, .
\end{multline*}
So, overall, 
\[
 \left \Vert  Y_{0,\ell}'  \sum_{i=0}^{k-1}\E_0 (Y_{i,\ell}') \right \Vert_{r/2}\leq 8  \left ( \int_0^1 (  \alpha_{2, {\bf Y}}^{-1} (u) \wedge k)^{r/2} Q^{r} (u \vee v) du \right )^{2/r}   \, .
\]
Taking into account this  bound in \eqref{TH25RIO} and using the fact that $\sum_{\ell =1}^L |a_{L} | \leq 1$, we derive
\[
\Vert  S'_q  \Vert^r_r \leq (16 r q)^{r/2} \int_0^1 (  \alpha_{2, {\bf Y}}^{-1} (u) \wedge q)^{r/2} Q^{r} (u \vee v) du\, .
\]
Since  $\alpha_{2, {\bf Y}}^{-1} (u) \wedge q = q \I_{0 <u \leq v} + (\alpha_{2, {\bf Y}}^{-1} (u) \wedge n )  \I_{v <u \leq 1}$, the above upper bound gives
\[
\frac{n}{q x^r} \Vert  S'_q  \Vert^r_r \ll  \frac{n}{ x^r}  ( qQ(v))^{r-1} \int_0^v Q(u) du +  \frac{n}{ x^r}  ( qM)^{r/2-1} \int_v^1 (  \alpha_{2, {\bf Y}}^{-1} (u) \wedge n)^{r/2} Q^{r/2 +1} (u ) du \, ,
\]
which proves \eqref{majorationSqr} by taking into account \eqref{restq}. 

We show now that 
\beq \label{majorationSq2}
 \frac{n^{r/2}}{q^{r/2}x^r} \Vert  S'_q  \Vert^r_2 \ll x^{-r}s^r_n(x) 
 +
 \left ( \frac n x \int_0^{L_n(x)} Q(u) du \right)^{r/2}
 \,  .
\eeq
By stationarity
\[
\E (   S'_{q}
)^2 =  \sum_{|i| \leq q } (q -|i|) \E (X_0' X_{|i|}') \, .
\]
Now, by Inequality (2.5) in \cite{DGM}, 
\[
\E (X_0' X_{|i|}') \leq 
4   \int_0^{\alpha_{1,{\bf Y}}(|i|) } Q^2 (u \vee v) du\, ,
\]
so that
\[
\E (   S'_{q})^2 \leq 8 q
\sum_{i=0}^{q-1} \int_0^{\alpha_{1,{\bf Y}}(i)}Q^2(u\vee v)du  
= 
8 q \int_0^1 (  \alpha_{1, {\bf Y}}^{-1} (u) \wedge q)
Q^2(u\vee v)du 
\, .
\]
Recall that $v=L_n(x)$ and that $qQ(v) \leq x$, by 
\eqref{restq}. Using in addition the fact that $\alpha_{1, {\bf Y}}^{-1} \leq \alpha_{2, {\bf Y}}^{-1} $, we get
\[
\E (   S'_{q})^2 \leq 8 q
 \int_{L_n(x)}^1 (  \alpha_{1, {\bf Y}}^{-1} (u) \wedge n)
Q^2(u)du 
+ 8 q x \int_0^{L_n(x)} Q(u) du
\, ,
\]
which proves \eqref{majorationSq2}.

Starting from \eqref{MarkovBigbloc} and taking into account \eqref{consRosenthal}, \eqref{majorationDkq}, \eqref{majorationSqr} and \eqref{majorationSq2}, we get 
\begin{multline}\label{ouf}
  \p \left (\max_{2 \leq 2j \leq [n/q]} \left | \sum_{i=1}^j \tilde U_{2i} \right | \geq x \right ) \ll \frac{s^r_n}{x^r}+   \frac{n}{x} \int_0^{L_n(x)}Q(u) du +  \frac{n}{x^{1+\beta/2}}   \int_0^{L_n(x)} R_n^{\beta/2}(u)   Q(u) du \\
  +   \frac{n}{x^{1+r/2}}   \int_{L_n(x)}^1 R_n^{r/2}(u)   Q(u) du \, .
  \end{multline}
Note that, to deal with  the second term on right hand in
 \eqref{majorationSq2}, we have used the fact if 
 $$ \frac{n}{x} \int_0^{L_n(x)}Q(u) du \geq 1$$ then the bound \eqref{ouf} is trivial, and otherwise 
 $$
\left ( \frac n x \int_0^{L_n(x)} Q(u) du \right)^{r/2} \leq 
  \frac n x \int_0^{L_n(x)} Q(u) du \, .
 $$
 
Obviously the upper bound \eqref{ouf} is also valid for the quantity
$$ 
\p \left (\max_{1 \leq 2j+1 \leq [n/q]} \left | \sum_{i=1}^j \tilde U_{2i+1} \right | \geq x \right ) \, .
$$
 Together with \eqref{dec15FN} and \eqref{decinter2} this completes the proof of the proposition for $q<n$.
 
The proposition is proved for any variable 
$X_i=f(Y_i)- {\mathbb E}(f(Y_i))$ with
$f=\sum_{\ell=1}^L a_{\ell} f_{\ell}$ and $f_{\ell} \in \tMon (Q, P_{Y_0})$,
$\sum|a_\ell|\leq 1$. 
Let us explain
how it can be extended to 
$\tMonm (Q, P_{Y_0})$. 

Let $f \in \tMonm(Q,P_{Y_0})$. By definition of
$\tMonm(Q,P_{Y_0})$, there exists $f_L= \sum_{\ell=1}^L a_{\ell, L}
g_{\ell, L}$ with $g_{\ell, L}$ belonging to $\tMon(Q, P_{Y_0})$ and
$\sum_{\ell=1}^L |a_{\ell, L}| \leq 1$, and such that $f_L$ converges
in ${\mathbb L}^1(P_{Y_0})$ to $f$. It follows that $X_{i,L}=f_L(Y_i)
- \E ( f_L(Y_i))$ converges in ${\mathbb L}^1$ to $X_i$ as $L$ tends
to infinity. Extracting a subsequence if necessary, one may also
assume that the convergence holds almost surely.
Hence, for any fixed $n$, $S_{n,L}=\sum_{k=1}^n X_{k,L}$ converges almost surely
and in ${\mathbb L}^1$ to $S_n$.

Let $Z_n= \max_{1 \leq k
\leq n} |S_{k}|$.  By Beppo-Levi,
\begin{equation}\label{BeppoLevy}
 \p \left ( \max_{1 \leq k \leq n} |S_k| > x
\right ) =
\lim_{k \rightarrow \infty} \p \left ( Z_n > x+k^{-1}
\right ).
\end{equation}
Let $h_k$ be a continuous function from ${\mathbb R}$ to $[0, 1]$, such that
$h_k(t)=1$ if $t >x+ k^{-1}$ and $h_k(t)=0$ if $t < x$. Let
$Z_{n,L}= \max_{1 \leq k \leq n} |S_{k, L}|$.
By Fatou's lemma,
\begin{equation}\label{Fatou}
\p \left ( Z_n > x+k^{-1}
\right )\leq \E \left (
 h_k(Z_n)
\right ) 
\leq \liminf_{L \rightarrow \infty}
\E \left (
 h_k(Z_{n, L})
\right )\leq \liminf_{L \rightarrow \infty}
\p \left (
Z_{n, L}>x
\right )\, .
\end{equation}
Now, we have just proved that  $\p \left (
Z_{n, L}>x
\right )$ satisfies the upper bound (\ref{eqmainmax}), uniformly 
in $L$.
From \eqref{BeppoLevy} and~\eqref{Fatou}, we infer
that $\p  ( \max_{1 \leq k \leq n} |S_k| > x
 ) $ satisfies also (\ref{eqmainmax}), which completes the proof 
 of Proposition \ref{deviationRosenthal}.

\section{Proofs of the results of Section \ref{Sec:new}}
\label{Sec:6}

\subsection{Preparatory material}

\begin{lma}\label{trivial} Let $p>1$.
 Assume that \eqref{WM} holds and let $b \in (0, p-1)$ and 
 $c>p-1$. Then
 $$
 \sup_{x>0}x^{p-1-b} \int_0^1 R^b(u) Q(u) {\bf 1}_{R(u)>x} du  < \infty \, , \quad 
 \text{and} \quad 
 \sup_{x>0}\frac 1 {x^{c-p+1}} \int_0^1 R^c(u) Q(u) {\bf 1}_{R(u)\leq x} du < \infty \, .
 $$
 \end{lma}

\begin{lma}\label{trivialbis} Let $p>1$.
 Assume that \eqref{WM0} holds and let $b \in (0, p-1)$ and 
 $c>p-1$. Then
 $$
 \lim_{x \rightarrow \infty}
 x^{p-1-b} \int_0^1 R^b(u) Q(u) {\bf 1}_{R(u)>x} du  =0 \, , \quad 
 \text{and} \quad 
 \lim_{x \rightarrow \infty} 
 \frac 1 {x^{c-p+1}} \int_0^1 R^c(u) Q(u) {\bf 1}_{R(u)\leq x} du =0 \, .
 $$
 \end{lma}
 \noindent{\bf Proof.}
Let $Z$ be a random variable with values in $(0,1)$, whose distribution 
has density
$$
u \rightarrow \bar Q (u) = \frac{Q(u)}{ \int_0^1 Q(u) du} \, .
$$
Then Condition \eqref{WM} is equivalent to 
\begin{equation}\label{WM2}
\sup_{x >0} x^{p-1}{\mathbb P}(R(Z)>x)  < \infty
\end{equation}
and Condition \eqref{WM0} is equivalent to 
\begin{equation}\label{WM02}
\lim_{x \rightarrow \infty} x^{p-1}{\mathbb P}(R(Z)>x) =0 \, .
\end{equation}
In the same way, the first statements of  Lemmas \ref{trivial} and  \ref{trivialbis} read respectively
\begin{equation}\label{stat1}
\sup_{x >0} x^{p-1-b}{\mathbb E}\left( (R(Z))^b{\bf 1}_{R(Z)>x}\right )  < \infty \, ,
\quad \text{and}
\quad 
\lim_{x \rightarrow \infty} x^{p-1-b} {\mathbb E}\left( (R(Z))^b{\bf 1}_{R(Z)>x}\right )=0 \, .
\end{equation}
Applying Fubini's theorem, one easily sees that
\begin{equation}\label{Fub1}
{\mathbb E}\left( (R(Z))^b{\bf 1}_{R(Z)>x}\right )= x^b {\mathbb P}(R(Z)>x) + 
b\int_x^\infty u^{b-1} {\mathbb P}(R(Z)>u) \, du \, .
\end{equation}
Now, \eqref{WM2} combined with \eqref{Fub1} implies the first part of \eqref{stat1},
and \eqref{WM02} combined with \eqref{Fub1} implies the second
 part of \eqref{stat1}. This concludes the proofs of the first statements of 
  Lemmas \ref{trivial} and  \ref{trivialbis}. 

To prove the second statements of  Lemmas \ref{trivial} and \ref{trivialbis}, we first note that they are respectively equivalent to 
\begin{equation}\label{stat2}
\sup_{x>0}\frac 1 {x^{c-p+1}} {\mathbb E} \left ( (R(Z))^c  {\bf 1}_{R(Z)\leq x} 
\right ) < \infty \, ,
\quad \text{and}
\quad 
\lim_{x \rightarrow \infty} 
\frac 1 {x^{c-p+1}} {\mathbb E} \left ( (R(Z))^c  {\bf 1}_{R(Z)\leq x} 
\right )=0 \, .
\end{equation}
Applying Fubini's theorem, one easily sees that
\begin{equation}\label{Fub2}
{\mathbb E}\left( (R(Z))^c{\bf 1}_{R(Z)\leq x}\right )= c \int_0^\infty 
 u^{c-1} {\mathbb P}(u<R(Z)\leq x) \, du \leq 
 c \int_0^x
 u^{c-1} {\mathbb P}(R(Z)>u)\, du \, .
\end{equation}
Now, \eqref{WM2} combined with \eqref{Fub2} implies the first part of \eqref{stat2},
and \eqref{WM02} combined with \eqref{Fub2} implies the second
part of \eqref{stat2}. This concludes the proofs of 
 of  Lemmas \ref{trivial} and  \ref{trivialbis}. 

\subsection{Proof of Theorem \ref{ineRosalphafaible}}

 To prove \eqref{RB}, it suffices to write 
\[
\left \Vert
 \max_{ 1 \leq k \leq n}  |
  S_k | 
\right \Vert_p^p
  = p \int_0^{\infty} x^{p-1} \p \left( \max_{ 1 \leq k \leq n}  |
  S_k | \geq x \right)\,   dx
   \leq s_n^p + p \int_{s_n}^{\infty} x^{p-1} \p \left( \max_{ 1 \leq k \leq n}  |
  S_k | \geq x \right)  \, dx \, ,
\]
and to apply Inequality \eqref{eqmainmax} with $r-2 < \beta < 2p-2 < r <2p$ to bound the last integral. The result  follows by applying Fubini's theorem. 

\subsection{Proof of Theorem \ref{Prop:LD}}

We apply Proposition \ref{deviationRosenthal} with 
$r/2=a \in (p-1,p)$ and $\beta \in (r-2, 2p-2)$. 

\noindent {\em Proof of Item 1.} 
We start by proving \eqref{WB}.
Recall that $s_n^2$ has been defined in \eqref{sn2}
and that $s_n^2(x)\leq s_n^2$. We have that
\begin{equation}\label{sn}
s^2_n \leq C(\alpha, Q) n\, , \quad  
\text{with} \quad 
C(\alpha, Q)=\sum_{i=0}^{\infty} \int_0^{\alpha_{1,
{\bf Y}} (i)}Q^2 (u) du \, .
\end{equation}
Since $p>2$, the condition \eqref{WM} implies that 
$
C(\alpha, Q)< \infty 
$. 
 Inequality  \eqref{eqmainmax}  and Remark \ref{rmk:sn2} yield
\begin{multline}\label{direct}
  \p \left( \frac 1 n  \max_{ 1 \leq k \leq n}  |
  S_k | \geq x \right)  \ll 
  \frac 1 {n^a x^{2a}} +   \frac{1}{x} \int_0^1Q(u)
   {\bf 1}_{R(u)>nx} du \\
    +  \frac{1}{n^{\beta/2} x^{1+\beta/2}}   \int_0^{1} R^{\beta/2}(u)   Q(u) {\bf 1}_{R(u)>nx} du \\
  +   \frac{1}{n^{r/2}x^{1+r/2}}   \int_{0}^1 R^{r/2}(u)   Q(u) {\bf 1}_{R(u)\leq nx} du  \, .
 \end{multline}
 From \eqref{WM}, we infer that
 $$
  \frac{1}{x} \int_0^1Q(u)
   {\bf 1}_{R(u)>nx} du \ll \frac 1 {n^{p-1} x^p} \, .
 $$
 To handle the two last terms on right hand in 
 \eqref{direct}, we apply Lemma \ref{trivial}. We then infer that
 \begin{align*}
 \frac{1}{n^{\beta/2} x^{1+\beta/2}}   \int_0^{1} R^{\beta/2}(u)   Q(u) {\bf 1}_{R(u)>nx} du 
 &\ll \frac 1 {n^{p-1} x^p} \, , \\
 \frac{1}{n^{r/2}x^{1+r/2}}   \int_{0}^1 R^{r/2}(u)   Q(u) {\bf 1}_{R(u) \leq nx} du &\ll \frac 1 {n^{p-1} x^p} \, ,
 \end{align*}
 and \eqref{WB} follows. 
 
 The proof of \eqref{WB2} is almost identical. 
 The only difference is about the first term on right hand in \eqref{direct}. This first term is equal to $(nx)^{-2a}s^{2a}_n(nx)$, 
  and for $c \in (0,1)$, one has
 $$
\frac{s^{2a}_n(nx)}{(nx)^{2a}} \leq \frac{1}{n^ax^{2a}}
\left( \int Q(u) R^c(u) R^{1-c}(u) {\mathbf 1}_{R(u) \leq nx} du 
\right )^a \, .
 $$
 Since $p=2$, it follows that $\int Q(u) R^c(u) du < \infty$, and finally
 $$
 \frac{s^{2a}_n(nx)}{(nx)^{2a}}  \ll \frac 1 {n^{ac} x^{a(1+c)}}\, ,
 $$
 completing the proof of \eqref{WB2}.

 \medskip
 
\noindent{\em Proof of Item 2.} We start from
\eqref{direct}.  Since the series with terms $n^{p-2}$ 
and $n^{p-2- (\beta/2)}$ are divergent, we obtain by 
summing in $n$,
$$
\sum_{n=1}^\infty  \frac{n^{p-2}}{x} \int_0^1Q(u)
   {\bf 1}_{R(u)>nx} du  \ll 
   \frac 1 { x^p} \int_0^1 (\alpha_{2,{\bf Y}}^{-1} (u))^{p-1}   Q^{p}(u) du  \, . 
$$
and 
$$
\sum_{n=1}^\infty  \frac{n^{p-2}}{n^{\beta/2} x^{1+\beta/2}}   \int_0^{1} R^{\beta/2}(u)   Q(u) {\bf 1}_{R(u)>nx} du   \ll 
   \frac 1 { x^p} \int_0^1 (\alpha_{2,{\bf Y}}^{-1} (u))^{p-1}   Q^{p}(u) du \, . 
$$
Since the series with term $n^{p-2- (r/2)}$ converges, 
we obtain by summing in $n$,
$$
\sum_{n=1}^\infty  
\frac{n^{p-2}}{n^{r/2}x^{1+r/2}}   \int_{0}^1 R^{r/2}(u)   Q(u) {\bf 1}_{R(u) \leq nx} du
  \ll 
   \frac 1 { x^p} \int_0^1 (\alpha_{2,{\bf Y}}^{-1} (u))^{p-1}   Q^{p}(u) du \, .
$$
The proof of Item 2 is now complete.

\subsection{Proof of Theorem \ref{Holder}}
Since $p>2$, 
the condition \eqref{WM0} implies \eqref{DMR}, and the process
$
\left \{ W_n(t), t \in [0,1]\right \}
$
converges in distribution in the space $C([0,1], \|\cdot \|_\infty)$ to  $\sigma W$. 
It remains to prove the tightness in ${\mathcal H}^0_{\delta}([0,1])$. We start from the tightness 
criterion given in \cite{Gir}: the process $
\left \{ W_n(t), t \in [0,1]\right \}
$ is tight in ${\mathcal H}^0_{\delta}([0,1])$ as soon as, 
for any $\varepsilon>0$,
\begin{equation}\label{Gir2}
\lim_{\eta \rightarrow 0} \limsup_{n \rightarrow \infty} 
n \sum_{k=1}^{[\log (n\eta)]} \frac 1 {2^k} {\mathbb P} 
\left ( \max_{1 \leq i \leq 2^k} |S_i|> \varepsilon 2^{k \delta} n^{1/p} \right) =0 \, ,
\end{equation}
where $\log$ is the binary logarithm. Applying Proposition 
\ref{deviationRosenthal}
with $s_n$ as in \eqref{sn}, $r \in (2p-2,2p)$, and  
$\beta \in (r-2, 2p-2)$, we get that
\begin{equation}\label{devbound}
n \sum_{k=1}^{[\log (n\eta)]} \frac 1 {2^k} {\mathbb P} 
\left ( \max_{1 \leq i \leq 2^k} |S_i|> \varepsilon 2^{k \delta} n^{1/p} \right)
\leq \sum_{i=1}^4 I_i(n, \eta) \, ,
\end{equation}
with 
\begin{align*}
I_1(n, \eta) &\ll \frac 1 {n^{(r-p)/p}} \sum_{k=1}^{[\log (n\eta)]}  2^{k(r-p)/p} \, ,\\
I_2(n, \eta) &\ll \sum_{k=1}^{\infty}  \frac {n^{(p-1)/p}} {2^{k \delta}}
\int_0^1 Q(u) {\bf 1}_{R(u)>\varepsilon 2^{k \delta} n^{1/p} } du \, ,\\
I_3(n, \eta) &\ll \sum_{k=1}^{\infty}  \frac {n^{(2p-2-\beta)/2p}} {2^{k \delta(2+ \beta)/2}}
\int_0^1 R^{\beta/2}(u)Q(u) {\bf 1}_{R(u)>\varepsilon 2^{k \delta} n^{1/p} } du \, ,\\
I_4(n, \eta) &\ll \sum_{k=1}^{\infty}  \frac {n^{(2p-2-r)/2p}} {2^{k \delta(2+ r)/2}}
\int_0^1 R^{r/2}(u)Q(u) {\bf 1}_{R(u)\leq \varepsilon 2^{k \delta} n^{1/p} } du
\, .
\end{align*}
Clearly $I_1(n, \eta) \ll \eta^{(r-p)/p}$ and consequently 
\begin{equation}\label{LI1}
\lim_{\eta \rightarrow 0} \limsup_{n \rightarrow \infty} I_1(n, \eta) =0 \, .
\end{equation}
Now 
$$
I_2(n, \eta) \ll  n^{(p-1)/p} 
\int_0^1 Q(u) {\bf 1}_{R(u)>\varepsilon n^{1/p} } du \, ,
$$
and it follows from \eqref{WM0} that
\begin{equation}\label{LI2}
 \lim_{n \rightarrow \infty} I_2(n, \eta) =0 \, .
\end{equation}
In the same way 
\begin{equation}\label{I3}
I_3(n, \eta) \ll   n^{(2p-2-\beta)/2p}
\int_0^1 R^{\beta/2}(u)Q(u) {\bf 1}_{R(u)>\varepsilon  n^{1/p} } du
\, .
\end{equation}
 Applying Lemma \ref{trivialbis} with $b=\beta/2$, it follows from \eqref{I3} that
 \begin{equation}\label{LI3}
 \lim_{n \rightarrow \infty} I_3(n, \eta) =0 \, .
\end{equation}
To control $I_4(n, \eta)$, we first note that, by Lemma \ref{trivialbis} with $c=r/2$,
$$
\frac {n^{(2p-2-r)/2p}} {2^{k \delta(2+ r)/2}}
\int_0^1 R^{r/2}(u)Q(u) {\bf 1}_{R(u) \leq \varepsilon 2^{k \delta} n^{1/p} } du
\ll \frac 1 {2^{k p \delta}} \, ,
$$ 
and 
$$
\lim_{n \rightarrow \infty} \frac {n^{(2p-2-r)/2p}} {2^{k \delta(2+ r)/2}}
\int_0^1 R^{r/2}(u)Q(u) {\bf 1}_{R(u)\leq \varepsilon 2^{k \delta} n^{1/p} } du =0 \, .
$$
Hence, by the dominated convergence theorem, 
\begin{equation}\label{LI4}
\lim_{n \rightarrow \infty} I_4(n, \eta) =0 \, .
\end{equation}

Finally, the tightness follows from \eqref{Gir2}, \eqref{devbound}, 
\eqref{LI1}, \eqref{LI2}, \eqref{LI3} and \eqref{LI4}.
This completes the proof of Theorem \ref{Holder}.

\end{document}